\documentclass{article}
\usepackage[totalwidth=380pt, totalheight=600pt]{geometry}
\usepackage{amsmath,amssymb,amsthm}
\usepackage[dvips]{graphicx,epsfig}

\usepackage[latin1]{inputenc}
\usepackage{amsmath,amssymb,amsthm}
\usepackage[dvips]{graphicx,epsfig}
\newcommand{\A}{\mathbb{A}}

\newcommand{\HH}{\mathbb{H}}

\newcommand{\NN}[0]{\mathbb{N}}
\newcommand{\PP}{\mathbb{P}}

\newcommand{\ZZ}{\mathbb{Z}}

\newtheorem{nt}{Notation}
\newtheorem{conjecture}[nt]{Conjecture}
\newtheorem{prop}[nt]{Proposition}

\newcommand{\f}{\frac}

\newcommand{\fd}{\ensuremath{\rightarrow}}

\newcommand{\findem}{\hfill\rule{2mm}{2mm}}

\renewcommand{\phi}{\ensuremath{\varphi}}

\newcommand{\inj}{\ensuremath{\hookrightarrow}}
\newcommand{\inc}{\ensuremath{\subset}}
\newcommand{\nl}{\ \\[2mm]}

\newcommand{\pla}{\A^2}
\newcommand{\plp}{\PP^2}
\renewcommand{\phi}{\ensuremath{\varphi}}

\newcommand{\s}{Spec\;}

\newcommand{\x}{\ensuremath{\times}}

\newcommand{\barr}{\begin{array}}
\newcommand{\earr}{\end{array}}
\newcommand{\bit}{\begin{itemize}}
\newcommand{\eit}{\end{itemize}}
\newcommand{\beq}{\begin{eqnarray*}}
\newcommand{\eeq}{\end{eqnarray*}}
\newcommand{\beqn}{\begin{eqnarray}}
\newcommand{\eeqn}{\end{eqnarray}}
\newcommand{\bconj}{\begin{conjecture}}
\newcommand{\econj}{\end{conjecture}}
\newcommand{\bcor}{\begin{coro}}
\newcommand{\ecor}{\end{coro} \noindent}
\newcommand{\ben}{\begin{enumerate}}
\newcommand{\een}{\end{enumerate} \noindent}
\newcommand{\bnot}{\begin{nt} }
\newcommand{\enot}{\end{nt} \noindent}
\newcommand{\bdefi}{\begin{defi}}
\newcommand{\edefi}{\end{defi} \noindent}
\newcommand{\bprop}{\begin{prop}} 
\newcommand{\brap}{\begin{rappel}}
\newcommand{\erap}{\end{rappel} \noindent }
\newcommand{\brq}{\begin{rem}}
\newcommand{\erq}{\end{rem} \noindent }
\newcommand{\bthm}{\begin{thm}}
\newcommand{\blm}{\begin{lm}}
\newcommand{\bex}{\begin{ex}}
\newcommand{\eex}{\end{ex}\noindent }
\newcommand{\bexo}{\begin{exo} \normalfont}
\newcommand{\eexo}{\end{exo}\noindent }

\font \tengothic=eufm10 scaled\magstep 1
\font\sevengothic=eufm7 scaled\magstep 1
\newfam\gothicfam
      \textfont\gothicfam=\tengothic
      \scriptfont\gothicfam=\sevengothic

\newtheorem{coro}[nt]{Corollary} 
\newtheorem{defi}[nt]{Definition}

\newtheorem{ex}[nt]{Example} 
\newtheorem{exo}[nt]{Exercise} 
\newtheorem{lm}[nt]{Lemma} 
\newtheorem{rappel}[nt]{Recall} 
\newtheorem{rem}[nt]{Remark} 
\newtheorem{thm}[nt]{Theorem}

\newcommand{\demo}{\noindent \textit{Proof}}    
\newcommand{\elm}{\end{lm} \noindent \textit{Proof: }}  
\newcommand{\eprop}{\end{prop} \noindent \textit{Proof: }} 
\newcommand{\ethm}{\end{thm}\noindent \textit{Proof: }}

\begin{document}
\sloppy
\title{Irreducible components of the equivariant punctual Hilbert schemes}
\date{June 10, 2001}
\author{Laurent Evain  (evain@tonton.univ-angers.fr) }
\bibliographystyle{plain}
\maketitle
\noindent
{\bf Abstract:}
Let $\HH_{ab}$ be the 
equivariant Hilbert  scheme  parametrizing 
the 0-dimensional subschemes of the affine plane invariant under the natural action 
of the one-dimensional torus $T_{ab}:=\{(t^{-b},t^a)\ t\in k^*\}$. 
We compute the irreducible components of $\HH_{ab}$: they are in
one-one correspondence with a set of Hilbert functions. 
\\
As a by-product of the proof, we give new proofs of results by
Ellingsrud and Str\o mme, namely the main lemma of the computation of the Betti numbers
of the Hilbert scheme $\HH^l$ parametrizing the 0-dimensional subschemes 
of the affine plane of length $l$
\cite{ellingsrud-stromme87:chow_group_of_hilbert_schemes}, 
and a description of 
Bialynicki-Birula cells on $\HH^l$ 
by means of explicit flat families 
\cite{ellingsrud-stromme88:cartes_sur_le_schema_de_Hilbert}.
In particular, we precise conditions of applications of this last description. 

\section*{Introduction}
\label{introduction}
The Hilbert
scheme $\HH$ parametrizing subschemes of the plane  $\s k[x,y]$ is a disjoint
union of its components $\HH^l$ parametrizing the subschemes of length
$l$. Two auxiliary subschemes are useful to study $\HH$. The first one
is $\HH_{*}$ parametrizing the subschemes supported by a point
$*$. For instance, the papers by Lehn and Nakajima  (\cite{lehn99:_chern_classes_of_tautological_sheaves_on_Hilbert_schemes} and 
\cite{nakajima97:_heisenberg_et_Hilbert_schemes}) 
are explicit situations where the geometry of $\HH_{*}$ 
has proven to be useful.
The second auxiliary scheme $\HH_{ab}$ encode in
some sense the action of the torus on $\HH$.  
The action we are talking about is the action 
of the two dimensional torus $k^*\x k^*$
induced by the linear action 
$(t_1,t_2).x^{\alpha}y^{\beta}=(t_1.x)^{\alpha}(t_2.y)^{\beta}$ of $k^*\x k^*$
on $k[x,y]$. If 
$T_{ab}=\{(t^{-b},t^a),\ t\in k^*\}$ is a one dimensional sub-torus,
the closed subscheme 
$\HH_{ab}\inc \HH$ parametrizes by definition the zero dimensional 
subschemes invariant under  
the action of $T_{ab}$. Its usefulness has been illustrated by 
Ellingsrud and Str\o mme's work 
\cite{ellingsrud-stromme87:chow_group_of_hilbert_schemes} where they compute 
the Chow group of $\HH^l$ by studying the equivariant embedding  
$\HH_{ab}\inj \HH$ for $(a,b)$ general enough.
See Brion \cite{brion97:_equivariant_chow_groups} 
for other links between $\HH$ and $\HH_{ab}$ relative to 
equivariant cohomology.
\\
So, even if
one is primarily interested in studying $\HH$, it is natural to study 
the three Hilbert schemes $\HH$, $\HH_{ab}$ and $\HH_*$ as a whole
for their respective geometries are linked.  
\nl
One of the basic questions about these schemes is to  determine their 
irreducible components. The answers are known for $\HH$ (the irreducible
components are the smooth subschemes $\HH^l$ by Hartshorne
\cite{hartshorne66:connexite_du_schema_de_hilbert} and Fogarty
\cite{fogarty68:lissite_du_hilbert} ) and for $\HH_*$ (the irreducible
components are $\HH_*^l:=\HH_* \cap \HH^l$ by Brian\c con
\cite{briancon77:description_de_Hilb_a2}),
but unknown for $\HH_{ab}$. 
\\
As to $\HH_{ab}$, it cannot be irreducible 
for there are obviously disconnected subschemes  $\HH_{ab}^l:=
\HH_{ab}\cap \HH^l$. These  can still be divided into smaller
disconnected pieces as follows. The ideal of a  
subscheme corresponding to a point 
of $\HH_{ab}^l$ is quasi-homogeneous (homogeneous if $a=1,b=-1$) 
and is a direct sum $I=\oplus_{l  \in \ZZ} I_l$ 
of vector spaces.  The codimensions $d_i$
of the $I_i$ are fixed 
on a connected component and verify $\sum d_i=l$. 
Thus $\HH_{ab}^l=\coprod H_{ab}(H)$, where 
$H:\ZZ\fd \NN$ runs through the functions verifying $\sum_{i\in \ZZ}
H(i)=l$ and 
$H_{ab}(H)$ parametrizes the zero-dimensional subschemes satisfying 
$(\dots,d_{-1},d_0,d_1,d_2,\dots)=H$. So the natural question concerning 
the irreducibility of the equivariant Hilbert schemes is:
are the natural candidates $\HH_{ab}(H)$  connected and irreducible ?

\begin{thm}
  \label{thm principal d'irreductibilite}
  For any $H$, if  $\HH_{ab}(H)\neq \emptyset$, then $\HH_{ab}(H)$ 
  is smooth and connected (hence irreducible).
\end{thm}
\noindent
This result was already known  for $(a,b)$ general enough
(since then $\HH^l_{ab}$
is a union of points) and  for the  special value
$(a,b)=(1,-1)$ thanks to Iarrobino \cite{iarrobino77:_punctual_hilbert_schemes_AMS}.
\nl
\textbf{On  the proof}.\\
The first point is that  $\HH_{ab}$ is smooth as it is the fixed locus 
of the smooth scheme $\HH$ under the action of the one-dimensional
torus $T_{ab}$
\cite{iversen72:lissite-de-partie-invariante-sous-l'action-d'un-tore}, 
so
we just have to prove the connectedness of $\HH_{ab}(H)$.
\\
The general plan is as follows. We explain that $\HH_{ab}(H)$
is naturally stratified by locally closed subschemes $C(E)$. The strata
are parametrized by staircases (or Young diagrams) $E$ and in each
stratum $C(E)$ there is a particular subscheme $Z(E)$. Each stratum is
connected as an affine space, so the problem is to link the strata 
together. This is done in two steps. Firstly, we partially order the strata and we 
show that each stratum is connected to 
a smaller stratum except for some special 
strata characterized by their staircase $E$. 
We then show that there is in fact a unique special stratum.
It follows that all strata can be connected to this
special stratum and the theorem follows.
\\ 
The main point is to find a workable condition to 
connect a stratum $C(E)$ to a smaller stratum.
To this purpose, we  describe the tangent space $T_{Z(E)}^{ab}$ 
to $\HH_{ab}(H)$ at the special point $Z(E)$ of the stratum.
It is a vector space whose a base is  
a combinatorial datum associated with $E$ 
(theorem \ref{thm:description de l'espace tangent invariant}). 
This vector space splits into a direct sum 
of a positive tangent space $T_{Z(E)}^{ab+}$ 
and a negative tangent space $T_{Z(E)}^{ab-}$. 
We describe a sort of  ``exponential'' map $T_{Z(E)}^{ab+}\fd \HH_{ab}$, which turns 
out to be an isomorphism on $C(E)$. 
Using this map, a point $p$ of $T_{Z(E)}^{ab+}$ corresponds
to a subscheme $X(p)$ of $\pla$ whose equations are completely described. 
We check that for a suitable action 
of $k^*$, if $p\neq 0$, then $lim_{t \fd 0} t.X(p)$ is a subscheme in a
smaller stratum (proposition 
\ref{prop:description de la degenerescence}). 
This allows us to connect a stratum $C(E)$ to a
smaller stratum provided the positive tangent space in $Z(E)$ is non trivial.
\\
To conclude, we exhibit an explicit staircase $E_m$ characterized
combinatorially by a property of minimality 
(theorem \ref{thm:description escalier final}). We show that 
going down from a stratum to a smaller stratum as explained above, 
the process always ends in $C(E_m)$. 
The existence of a minimal $E_m$
is a purely
combinatorial statement. Its proof, 
however, requires the
characterization of minimality using the positive
tangent space.
The method here is very much in the spirit of toric
varieties \cite{fulton93:varietes-toriques}
where combinatorial statements are proved via
algebro-geometric arguments.
\nl
\textbf{By-products of the proof.}
The methods involved in the proof to study $\HH_{ab}$  remain useful
with few changes to study $\HH$. This allows us to recover results
by Ellingsrud-Str\o mme
(\cite{ellingsrud-stromme87:chow_group_of_hilbert_schemes} 
and \cite{ellingsrud-stromme88:cartes_sur_le_schema_de_Hilbert}). 
\\
More precisely, the
tangent space $T_{Z(E)}$ to $\HH$ at $Z(E)$ is a direct sum of the 
various invariant tangent spaces $T_{Z(E)}^{ab}$. 
Following the lines of \cite{ellingsrud-stromme87:chow_group_of_hilbert_schemes}, 
the control of $T_{Z(E)}$ is the key to compute
the Betti numbers of $\HH^l$ and we get by our description of $T_{Z(E)}$ 
a new proof of their main technical lemma. 
\nl
This tangent space $T_{Z(E)}$ splits into a direct sum $T_{Z(E)}^+\oplus T_{Z(E)}^-$,
just like the invariant tangent space $T_{Z(E)}^{ab}$. 
And we can still define an
embedding $e:T_{Z(E)}^+ \fd \HH$ using the same methods as for the
invariant case, that is by exhibiting an explicit flat family on $T_{Z(E)}^+$. 
The image of $e$ is a Bialynicki-Birula cell $C(E)$
with respect to a suitable action of the torus
$k^*$ on $\HH$ (theorem \ref{prop:e_est_un_plongement cas general}). 
In other words, we have explicit charts for cells $C(E)$.
Though it has been obtained by quite different methods, 
the flat family which defines $e$ is up to an
identification the same family as in
\cite{ellingsrud-stromme88:cartes_sur_le_schema_de_Hilbert}.
We note that the families involved in 
\cite{ellingsrud-stromme88:cartes_sur_le_schema_de_Hilbert}
are more general than ours since they give explicit charts for 
cells $C(E)$ defined by a more general action of the
torus $k^*$. However, some of the extra cases may have
gaps. Consequently, theorem \ref{prop:e_est_un_plongement cas general} 
can be seen as a prolongation of the analysis of Ellingsrud and
Str\o mme: it provides conditions 
on the action of $k^*$ under which 
the description of the cells given in 
\cite{ellingsrud-stromme88:cartes_sur_le_schema_de_Hilbert}
are valid (section \ref{ssec:exponentielle cas general}).

\section{The stratified subscheme $\HH_{ab}(H)$.} \label{sec:def de la stratification}
\subsection{Disconnected  subschemes of $\HH_{ab}$} \label{ssnot}
In this section, we introduce the closed subschemes $\HH_{ab}(H)\inc
\HH_{ab}$.
As claimed in the introduction, for $H_1\neq H_2$, 
$\HH_{ab}(H_1) \cap \HH_{ab}(H_2)=\emptyset$. 
\nl
Recall that 
$k$ is an algebraically closed field, that 
$a$ and $b$ are two 
relatively prime integers, 
and that $T_{ab}=\{(t^{-b},t^a), t \in k^*\}$ is a one 
dimensional subtorus of $k^*\x k^*=:T$. 
Since $(a,b)\neq (0,0)$  we can suppose by symmetry $b\neq 0$, and
since $(a,b)$ is defined up to sign, we suppose from now on $b<0$.
The Hilbert scheme $\HH^l$ 
parametrizing the zero-dimensional subschemes 
of $\s k[x,y]$ of length $l$ is connected and by definition 
$\HH=\amalg\; \HH^l $ 
where $\amalg$ stands for a disjoint  union.
There are actions  $T\x \HH \fd \HH$ and $T_{ab}\x \HH \fd \HH$ 
of the tori on  
$\HH$, which in terms of coordinates are given  by  
$(t_1,t_2).(x^{\alpha}y^{\beta})=(t_1x)^{\alpha}(t_2y)^{\beta}$.
We denote by $\HH_{ab}\subset \HH$  the subscheme  parametrizing
the subschemes of the plane invariant under the action of
$T_{ab}$. Then $\HH_{ab}=\amalg\; \HH_{ab}^l$, where 
$\HH_{ab}^l:=\HH_{ab} \cap \HH^l$. One can still separate $\HH_{ab}^l$ into disjoint 
subschemes by 
fixing a Hilbert function. 
\nl 
To do this, let us  characterize the subschemes of $\s k[x,y]$ which
lie in $\HH_{ab}$ by their ideals. 
Define the degree $d$ of a monomial by the formula
$d(x^\alpha y^\beta)=-b \alpha+ a \beta$. 
If $I$ is an ideal of $k[x,y]$, we let
$I_n:=I \cap k[x,y]_n$, where $k[x,y]_n$ denotes the vector space
generated by the monomials $m$ of degree $d(m)=n$. 
A subscheme $Z$  is in $\HH_{ab}$, 
if and only if its ideal  is quasi-homogeneous with respect to $d$, 
ie. $I(Z)=\oplus_{n \in \ZZ} I(Z)_n$.
\nl
By semi-continuity, if a subscheme $Z' \in \HH_{ab}$ is a 
specialization of $Z\in  \HH_{ab}$, then the codimensions 
of $I_n(Z)$ and $I_n(Z')$ in $k[x,y]_n$ verify 
$codim I_n(Z) \geq codim I_n(Z')$. 
But $Z$ and $Z'$, being in the same connected component of $\HH$,
have the same length
$l=\sum _{n}codim I_n(Z)= \sum _{n}codim I_n(Z')$.
It follows that the 
sequence $H(Z)=(\dots,h_{-1},h_0,h_1,\dots )$,
where $h_i=codim I_n(Z)$, is constant 
on the connected components 
of $\HH_{ab}$ and that $h_n=0$ for $n>>0$ and $n<<0$. 
If $H=(\dots,0,0,h_r, \dots,h_s,0,0,\dots )$ is any sequence, 
we denote by $\HH_{ab}(H)$ 
the closed subscheme of $\HH_{ab}$ (possibly empty) 
parametrizing the subschemes 
$Z$ verifying $H(Z)=H$. By the above, we have 
$\HH_{ab}=\amalg\;\HH_{ab}(H)$.

\subsection{The case $ab\geq 0$ via a theorem of Bialynicki-Birula} \label{le cas ab>=0}
In this section, we introduce the stratification on $\HH_{ab}(H)$,
the  special points $Z(E)$ of the cells $C(E)$ and we prove the theorem
under the condition
$ab\geq 0$.
The strata can 
be defined in terms of Gr\"obner bases, using Grassmannians 
or using a Bialynicki-Birula theorem, as explained in
\cite{evain00:cnincidence_cellules_schubert}.
We recall here the first and the last approach. The theorem follows
easily if $ab\geq 0$: 
the stratification on $\HH_{ab}(H)$ is reduced to 
one stratum, which is an affine space by Bialynicki-Birula.
\nl
Each cell is associated with a staircase $E$: 
to a staircase $E$ corresponds a subscheme $Z(E)$ and the cell
$C(E)$ is the unique cell containing the subscheme $Z(E)$. 
To be more precise, we recall that
a staircase is a subset of $\NN^2$ whose complementary
is stable by addition of $\NN^2$. In this paper, we will 
identify freely the monomial $x^p y^q$ with the couple $(p,q)$ 
and therefore the expression ``staircase of monomials'' will make sense.
More generally, we will transpose unscrupulously the definitions
between couples of integers and monomials. If $E$ is a staircase,
then the vector space  $I ^E$ generated by the monomials which are not 
in $E$ is an ideal and conversely, every monomial ideal is an 
ideal $I^E$ for a unique staircase $E$. The subscheme  $Z(E)$
whose ideal is $I ^E$ is in $\HH_{ab}(H)$ if and only if 
$E$ has $H(i)$ elements in degree $i$.
\nl
{\bf The Gr\"obner bases point of view.} 
We choose to order the monomials of $k[x,y]$ by the rule:
$m_1<m_2$ if $(d(m_1),d_y(m_1))<(d(m_2),d_y(m_2))$ 
for the lexicographic order, where 
$d_y$ is defined by
$d_y(x^{\alpha}y^{\beta})=\beta$.
 Let $m_1, m_2, \dots$ be the monomials which don't belong to $E$.  
Fixing once and for all $a,b$ and
$H=(\dots,0,0,h_r,\dots,h_s,0,\dots)$,
an ideal of $\HH_{ab}(H)$  
is in $C(E)$ if, regarding it as a $k$-vector space, it admits a basis  
$f_1, f_2, \dots$ where $f_i=m_i+R_i$, $R_i$ being a linear combination  
of monomials strictly smaller than $m_i$.   The locus $C(E)$ in
$\HH_{ab}$ is non
empty exactly when $E$ has $h_i$ elements in degree $i$.
\\
If $ab<0$, the above order is a monomial order in the sense 
of Gr\"obner bases. In this case,
the theory  
of Gr\"obner bases associates with every  
ideal in $k[x,y]$ a monomial ideal called initial ideal 
and $C(E)$ is the locus in $\HH_{ab}(H)$ parametrizing the  
ideals whose initial ideal is $I^E$.
\nl
{\bf The Bialynicki-Birula point of view.}
Let $X$ be a smooth projective variety 
over $k$
admitting an action of 
the torus $k^*$. Suppose that the action has a 
finite number of fixed points $x_1, \dots,x_n$. Let 
$T_{X,x_i}^+$ be the part of the tangent space to $x_i$ in 
$X$ where the weights of the $k^*$-action are positive, and 
let $X_i:=\{x \in X, \lim_{t \fd 0}(t.x)=x_i\}$. Then 
a theorem of Bialynicki-Birula asserts that the $X_i$ are
a cellular decomposition of $X$ in 
affine spaces and satisfy $T_{X_i,x_i}=T_{X,x_i}^+$.
If $U\inc X$ is a stable open subset of $X$, it is still possible 
to define cells which are affine spaces associated with the fixed
points which lie in $U$, but these cells do not always cover $U$.
In our case, fixing two integers 
$p$ and $q$ with $ap+bq>0$, 
the torus $k^*$ acts on $k[x,y]$ by $t.x=t^px$ and $t.y=t^{q}y$. 
This action induces an action of $k^*$ on $\HH_{ab}(H)$. 
The fixed points of $\HH_{ab}(H)$ under $k^*$ are 
the monomial subschemes $Z(E)$.
Applying the Bialynicki-Birula theorem to the action of $k^*$ on $\HH_{ab}(H)$,
we get a set of cells. We denote by $C(E)$ the cell associated with the
fixed  point $Z(E)$. The previous description of these cells 
insure that they cover $\HH_{ab}(H)$.
\nl
\textbf{Proof of the theorem in the case $ab\geq 0$.}
When the product $ab$ is non negative, there is at most one staircase
$E$ compatible  
with the Hilbert function $H$ (i.e. such that $E$ has $h_i$ elements in 
degree $i$). It follows that $\HH_{ab}(H)$ is empty or an affine
space $C(E)$.
The theorem is then obvious.

\section{Description of the tangent space}
Let $E$ be a staircase,  
$T_{Z(E)}$ be the tangent space to  $\HH$ at the point $Z(E)$,
$T_{Z(E)}^{ab}$ be the tangent space to $\HH_{ab}$ at the point $Z(E)$.
In this section, we give a description of $T_{Z(E)}^{ab}$ and
$T_{Z(E)}$.
\nl
We need some combinatorial vocabulary that we introduce now.
\\
\textbf{Clefts and cleft couples.} 
A cleft for $E$ is a couple $(u,v)$ such that $x^uy^v$
is a monomial in $I^E$ minimal for the divisibility relation among
the monomials of $I^E$. A positive (resp. negative)
half-direction is a couple of relatively prime integers $(f,g)$ with $f>0$
(resp. $f<0$) or $f=0$ and $g<0$ (resp. $g>0$). A couple
of points $(M,N)$ in $\NN^2$ has half-direction $(f,g)$ if the vector
$\vec{MN}$ is a positive multiple of $(f,g)$. 
We have a notion of direction by identifying  two opposites  half
directions. A cleft couple (relatively to $E$) with half-direction $h$
is a couple of elements $(c,m)$ in $\NN^2$ s.t. $c$ is a cleft, $m \in E$,
and $(c,m)$ has half-direction $h$.
\nl
\textbf{The orders $>_+$ and $>_-$ on cleft couples.}
We put 
orders on the set of clefts and on the set of cleft 
couples. The identification between $x^py^q$  and $(p,q)$
gives a lexicographic order $>_+$ on the monomials of $k[x,y]$. 
The identification between $x^py^q$  and $(q,p)$ gives a lexicographic order $>_-$.
Any order on the monomials induce an order on the clefts (by
restriction), on the couples of monomials (lexicographically) and 
on the cleft couples (the restriction of the latter). 
We still note the induced orders $>_+$ and $>_-$. 
\nl
\textbf{Some particular sets of cleft couples and associated vector
  spaces.}
The combinatorial object used to describe $T_{Z(E)}$ is the
significant cleft couple:
\begin{defi}
  Let $(c,m)$ be a cleft couple with positive 
(resp. negative) half-direction, $c'$ be the cleft successor
of $c$ for  $>_+$ (resp. for  $>_-$) and 
$s$ be  the smaller common multiple of $c$ and $c'$. 
Then $(c,m)$ is said to be 
a significant couple if $m\f{s}{c} \notin E$.
\end{defi}
\noindent
For a cleft couple, we'll be interested to know whether it is 
significant and what is its direction. Formally, we introduce the sets 
 $C,\overline C, C_{ab}, \overline {C}_{ab}$ containing respectively
\begin{itemize}
\item the significant cleft couples (with respect to $E$)
\item the cleft couples
\item the significant cleft couples with direction $(a,b)$.
\item the cleft couples with direction $(a,b)$
\end{itemize}
We have the obvious inclusions:
$$
\begin{array}{ccc}
C & \fd & \overline C\\
\uparrow & & \uparrow \\
C_{ab} & \fd & \overline{C}_{ab}
\end{array}
$$
For each of the above sets, we can form the vector space on this 
set, i.e. the vector space whose elements are the formal linear
combinations of elements of this set.
We denote by $R,\overline {R},R_{ab},\overline {R}_{ab}$ the vector spaces 
corresponding to  $C,\overline {C},C_{ab},\overline {C}_{ab}$
To the inclusion of sets corresponds the inclusion of vector spaces
$$
\begin{array}{ccc}
R & \fd & \overline R\\
\uparrow & & \uparrow \\
R_{ab} & \fd & \overline{R}_{ab}
\end{array}
$$
\nl
We can now formulate the main result of the section, whose proof ends 
with corollary \ref{coro description espace tangent}.
\begin{thm} \label{thm:description de l'espace tangent invariant}
There is a natural isomorphism $T_{Z(E)}^{ab}\simeq R_{ab}$. 
\end{thm}
We start by constructing an injective morphism
$\phi:T_{Z(E)}^{ab}\fd \overline R$. Next, we will identify the
image $\phi(T_{Z(E)}^{ab})$ with $R_{ab}$.
\nl
We have the classical description of tangent spaces on Hilbert 
schemes $T_{Z(E)}\simeq Hom_{k[x,y]}(I^E,k[x,y]/I^E)$. 
Let's recall in our situation how  
to produce the infinitesimal 
deformation of $k[x,y]/I^E$ starting from an element $f \in
Hom_{k[x,y]}(I^E,k[x,y]/I^E)$.
\\
Let $V_E$ be the vector space generated by the monomials 
which are in $E$, $\pi:V_E \fd k[x,y]/I^E$ the isomorphism 
induced by the restriction of the projection $k[x,y]\fd k[x,y]/I^E$, and 
$p:k[x,y]/I^E \fd k[x,y]$ the inverse of $\pi$. 
Let  $\overline{f}=p\circ f: I^E\fd k[x,y]$. 
The set of elements $J=\{m+\epsilon \overline{f}(m), m \in I^E\}$ 
is an ideal of $k[x,y][\epsilon]/(\epsilon^2)$.
The quotient
$k[x,y][\epsilon]/J$ is the flat infinitesimal deformation of
$k[x,y]/I^E$ corresponding to $f$.
\\
For every cleft couple $(c,m)$, define $\lambda_{c,m}$ by the 
formula $\overline{f}(c)=\sum \lambda_{c,m}m$.
\begin{prop}
The linear map
$$
\begin{array}{cccc}
\phi:& T_{Z(E)} &\fd & \overline{R}\\
     & f \mbox{ s.t. }\overline{f}(c)=\sum \lambda_{c,m}m & \mapsto & \sum \lambda_{c,m}(c,m)
\end{array}
$$
is injective and its restriction to 
$T_{Z(E)}^{ab}$ factorizes:
$$
\begin{array}[]{rlc}
T_{Z(E)}^{ab} & \fd & \overline{R}\\
& \searrow & \uparrow\\
& & \overline{R}_{ab}
\end{array}
$$
\end{prop}
\demo: the ideal 
$I^E$ being generated by the clefts, the morphism 
$f \in Hom_{k[x,y]}(I^E,k[x,y]/I^E)$ 
is characterized by the images of the clefts. But these 
images are themselves characterized by 
$\phi(f)$ so $\phi$ is injective. 
\\
 By
Bialynicki-Birula \cite{bialynicki-birula73:theorems_on_actions_of_groups_annals},
$T_{Z(E)}^{ab} \inc T_{Z(E)}$ is the subspace 
containing the vectors fixed under the action of $T_{ab}$.
If $f$ is in 
$T_{Z(E)}^{ab} \inc T_{Z(E)}$, one then sees from the description 
of $f$ as an infinitesimal deformation of $k[x,y]/I$ 
that $\lambda_{c,m}$ is different from zero only 
if $(c,m)$ has direction $(a,b)$. It follows 
that $\phi(T_{Z(E)}^{ab}) \inc \overline{R}_{ab}$.
\findem \\
The next step is to show 
the equality $\phi(T_{Z(E)}^{ab})=R_{ab}\inc \overline{R}_{ab}$.
This is done in two steps. Firstly, we construct a graph $G$ associated 
to $E$, a vector space $R_G$ from $G$ and we show the
equality $\phi(T_{Z(E)}^{ab})=R_G$. 
We then conclude with the isomorphism $R_G\simeq R_{ab}$.
\nl
We construct a graph $G$ from the staircase $E$. The set 
of points is the set $C_{ab}$ of cleft couples with direction 
$(a,b)$. Significant couples are not the end of any arrow.
Non significant cleft couples are the end of exactly one arrow.
Let $(c_1,m_1)$ be  non significant with positive half
direction, and $c_2$ be the successor cleft of $c_1$ for $>_+$. If
$m_2:=m_1.\f{c_2}{c_1} \in k[x,y]$, draw an arrow from $(c_2,m_2)$ to
$(c_1,m_1)$, otherwise draw an arrow from $(c_1,m_1)$ to itself.
Replacing $>_+$ by $>_-$, you get arrows ending at cleft couples 
with negative half-direction.  
\\ 
If $G$ is a graph with set of points ${\cal P}=
\{p_1, \dots, p_n\}$
and arrows $f_1,\dots,f_p$, denote $o(f_i)$ the origin 
of $f_i$ and $e(f_i)$ the end of $f_i$.  Let $R_{\cal P}$ 
be the vector space whose base is the set of elements of 
$\cal P$ and $R_G$ be the sub-vector space of $R_{\cal P}$ 
whose elements are the linear combinations  
$\sum \lambda_{p_i} p_i$ verifying 
$\lambda_{p_i}=\lambda_{p_j}$ if there is an arrow from $p_i$ to $p_j$
and $\lambda_{p_i}=0$ if there is an arrow from $p_i$ to itself.

\begin{prop}
If $a.b\leq 0$, 
then  $\phi(T_{Z(E)}^{ab})=R(G)$.
If $a.b>0$, 
then $\phi(T_{Z(E)}^{ab})$ is a subspace of $R(G)$.
\end{prop}
\demo: we have seen that an element $f$ in $T_{Z(E)}=Hom_{k[x,y]}(I^E,k[x,y]/I^E)$ 
is characterized by the images $g_i=f(c_i)$ of the clefts $c_i$. 
Reciprocally, if we prescribe an image $g_i$  for each
cleft $c_i$, a compatibility  relation insures the existence of  $f$ in $T_{Z(E)}$
sending $c_i$ to $g_i$: if $c_i<_+c_j$ are two clefts such that $c_j$ 
is the successor of $c_i$,  and
if $s$ is the smaller common multiple  (s.c.m.) of $c_i$ and $c_j$, then their
images $g_i$,$g_j$ have to verify 
$$
g_i.\frac{s}{c_i}=g_j.\frac{s}{c_j}
$$
The element $g_i$ is a linear combination $\sum \lambda_{i,k } \overline{m_k}$ 
and we identify it with the element $\sum \lambda_{i,k} (c_i,m_k)$ 
using $\phi$. We do similarly for $g_j$. The above relation
between $g_i$ and $g_j$ translates into  relations between 
the coefficients $\lambda_{i,k}$ and $\lambda_{j,l}$. 
One sees that these relations are the relations 
given by the arrows of the graph in the case $ab\leq 0$ so 
we are done. In the case $ab>0$, there are at least these 
relations and possibly some others. 
\findem

\begin{coro}
\label{coro description espace tangent}
If $a.b\leq 0$, then 
$T_{Z(E)}^{ab}\simeq R_{ab}$.
If  $a.b>0$, 
then $T_{Z(E)}^{ab}=R_{ab}=0$.
\end{coro}
\demo: by construction, each point in the graph $G$ is the end of at most 
one arrow and the origin of at most one arrow. It follows that the
connected components of the graph  $G$ are chains of consecutive
points $p_1,p_2,\dots,p_n$ where two points $p_i$ and $p_{i+1}$ are 
connected by one arrow from $p_i$ to $p_{i+1}$. Moreover the $p_i$ 
are distinct if $n\geq 2$. The
vector space $R(G)$ is then obviously isomorphic to
$R(G')$, where  $G'$ is the graph with no arrow 
obtained from $G$ by keeping the points that are not the end 
of any arrow. The set of points of $G'$ is 
just the set $C_{ab}$ of significant cleft couples 
with direction $(a,b)$ so $R(G)=R(G')=R_{ab}$ and the first 
claim of the corollary is a consequence of the previous proposition. 
In the case $a.b>0$, $G'$ is empty because there are no 
significant cleft couples. So $R(G)=0$ and the second 
claim follows again from the last proposition.
\findem 
\\
This corollary obviously concludes the proof of theorem 
\ref{thm:description de l'espace tangent invariant}.
It also clarifies the structure of $\HH_{ab}$ in some cases. 
We have already explained that $\HH_{ab}$ is a disjoint union
of affine spaces if $ab\geq 0$. These affine spaces are particularly
simple if $ab>0$: 
\begin{coro} \label{si ab>0 hilbert equiv=ens de points}
  If $ab>0$, then $\HH_{ab}(H)$ is empty or reduced to a point. 
\end{coro}
\demo: there is at most one staircase $E$ compatible with $H$. If $E$
does not exist, then $\HH_{ab}(H)=\emptyset$. Otherwise,
$\HH_{ab}(H)=C(E)$. The cell $C(E)$ is an affine space by
Bialynicki-Birula and its dimension is zero by the previous corollary 
(\ref{coro description espace tangent}) \findem
\\
Following the lines of the proof of theorem 
\ref{thm:description de l'espace tangent invariant} and simply 
forgetting the arguments concerning the action of the torus $T_{ab}$, 
we have a description of the tangent space $T_{Z(E)}$ instead of a
description of the invariant tangent space $T_{Z(E)}^{ab}$.
\begin{thm} \label{thm:description de l'espace tangent}
The tangent space $T_{Z(E)}$ to $\HH$ at $Z(E)$ is  isomorphic to the vector space 
$R$ whose base is the set of significant cleft couples of $E$. 
\end{thm}

\begin{rem}
  This description of the tangent space implies the main lemma (3.2)
  of \cite{ellingsrud-stromme87:chow_group_of_hilbert_schemes}
\end{rem}

\subsection{Definition of the positive and negative tangent spaces}
\label{sec:def_espace_tangent positif et negatif}
Let $C_{ab+}$ (resp. $C_{ab-}$) be the set of significant cleft
couples with direction $(a,b)$ and positive (resp. negative)
half-direction. Let $R_{ab+}$ (resp. $R_{ab-}$) be the vector space 
on the set  $C_{ab+}$ (resp. on $C_{ab-}$). The decomposition 
$C_{ab}=C_{ab+} \cup C_{ab-}$ gives $R_{ab+}\oplus R_{ab-}=R_{ab}$.
Considering the isomorphism $T_{Z(E)}^{ab}=R_{ab}=R_{ab+}\oplus
R_{ab-}$, we will say that $R_{ab+}$ (resp. $R_{ab-}$)  is the
positive 
(resp. negative) tangent space to
$\HH_{ab}$ at $Z(E)$ and we will write it down $T_{Z(E)}^{ab+}$
(resp. $T_{Z(E)}^{ab-}$). We define $T_{Z(E)}^{+}:=\oplus T_{Z(E)}^{ab+}$
where the sum runs over all the directions $(a,b)$.

\section{The exponential map}
\label{sec:l'application exponentielle}

\subsection{The invariant case}
\label{ssec:exponentielle cas invariant}
By the previous section, there is an isomorphism  
$T_{Z(E)}^{ab+} \simeq \s k[X_{c,m}]$ where the $X_{c,m}$ are variables in
bijection with the 
significant cleft couples $(c,m)$ with direction  $(a,b)$ and positive
half-direction. The goal of this section is to produce a sort of  
``exponential''  map $e:T_{Z(E)}^{ab+} \fd \HH_{ab}$ which induces an
isomorphism between $T_{Z(E)}^{ab+}$ and the image $C(E)$. In other
words, this section provides an explicit chart for the cell $C(E)$.
In the remainder of this section, we will make the assumption
that the half-direction $(a,b)$ verifies  $a>0$ and $b<0$.
\nl 
The morphism $T_{Z(E)}^{ab+} \fd \HH$ 
corresponds to a universal ideal over $\s k[X_{c,m}]$
that we describe now. 
Let $c_1<_+c_2<_+ \dots <_+ c_n$ be  the
clefts of $E$.  We define 
a set of monomials $\Delta_i$ for $1 \leq i\leq n$ by decreasing
induction: $\Delta_n$ is the set a monomials divisible by $c_n$ and, for
$i\neq n$,  $\Delta_i$ is the set of monomials divisible by $c_i$ but not
divisible  by $c_{i+1}$. We denote by $E_i$ the set of significant cleft couples $(c,m)$ 
verifying $c=c_i$. We define a polynomial $P(c_i)\in k[X_{c,m}][x,y]$ for each
cleft $c_i$ and a polynomial $Q(c_i,m)\in k[X_{c,m}][x,y]$ for  each
$(c_i,m)\in E_i$ by decreasing induction on $i$. When
$i=n$, we put  
$$P(c_n)=c_n.$$ 
The set $E_n$ being empty, there is no polynomial
$Q(c_n,m)$ to define.  For a general $i$ and a significant cleft couple
$(c_i,m)$,  $g$ the
s.c.m. of $c_i$ and $c_{i+1}$, let $c_{i+k} \
(k>0)$ be the cleft such that $m.(g/c_i)$ is in $\Delta_i$. We put
$$Q(c_i,m)=P(c_{i+k}).\f{m}{c_{i+k}}$$ and 
$$P(c_i)=P(c_{i+1}).\f{c_i}{c_{i+1}}+\sum_{m \mbox{ s.t. $(c_i,m)$ is
significant}} X_{c_i,m}Q(c_i,m).$$

\begin{prop}
\label{prop:platitude_de_la_famille}
a) The ideal $I_+$ generated by the $P(c_i)$ defines a flat family $Z_+$
of constant length over
$\s k[X_{c,m}]$. \\
b) The fiber of $Z_+$ over the origin is $Z(E)$.
\end{prop}
\demo: when $X_{c_i,m}=0$, for all $(c_i,m)$, then $P(c_i)=c_i$ and
$I_+=I^E$ 
This shows point b). \\
To verify a) we use the theory of Gr\"obner bases. It suffices to find a
monomial order on $Spec\  k[x,y]$ such that, with respect to that
monomial order, the $P(c_i)$ are a Gr\"obner basis
of $I_+$ over each point and verify $in(P(c_i))=c_i$. Indeed, if it is the 
case, the monomials which  are in $E$ form a basis of $k[x,y]/I_+$ 
over each point,
so the
length of the fibers of $k[x,y]/I_+$ equals the cardinal of
$E$. Each fiber being supported by the origin of $\pla$, the assertion
on the length of the fibers implies the flatness. 
We choose the monomial order $>_-$. 
According  to the Buchberger  algorithm, to verify that the $P_i$ 
form a Gr\"obner basis, we
must verify that if $c_i$ and $c_{i+1}$ are two  consecutive cleft couples,
and if $g$ is the s.c.m. of $c_i$ and $c_{i+1}$, 
then the remainder of a division of $h=\f{g}{c_i} P(c_i)- \f{g}{c_{i+1}}
P(c_{i+1})$ by the $P(c_i)$ is zero. But 
$h=\f{g}{c_i} \sum_{m \mbox{ s.t.
$(c_i,m)$ is significant}} X_{c_i,m}Q(c_i,m)$. By construction, 
$\f{g}{c_i}Q(c_i,m)$ is a product of a polynomial $P(c_{i+k})$ by 
a monomial. So the remainder of a division is zero.\findem

\begin{thm}
The morphism $e:T_{Z(E)}^{ab+}\fd \HH$ defined by $I_+$
induces an isomorphism $T_{Z(E)}^{ab+}\simeq C(E)$.
\label{prop:e_est_un_plongement} 
\end{thm}
\demo: by construction, the image of $e$ verifies  $Im(e)\inc C(E)$.
If $e$ is an embedding, then $Im(e)=C(E)$ and the theorem is proved.
Indeed, suppose that $e$ is an embedding and that there exists  $p \in
C(E)-Im(e)$. Then the orbit $T_{ab}.p$ is such that $T_{ab}.p \cup Im(e)$
is singular at $Z(E)$. Counting the dimensions of the tangent spaces
gives the inequality 
$$
dim T_{C(E),Z(E)} \geq \dim T_{T_{ab}.p \cup Im(e),Z(E)}>\dim
  T_{Z(E)}^{ab+}.
$$
This contradicts the equality of the first and last tangent spaces
asserted by Bialynicki-Birula.\\
We will conclude by showing using Pl\"ucker coordinates 
that $i\circ e: T_{Z(E)}^{ab+}\fd \PP^N$ 
is an embedding, where $i: \HH^{card(E)} \fd \PP^N$ is an embedding 
of the Hilbert scheme in a projective space. 
\\
Fix $n>0$ 
and define the set of monomials $C_n$ by 
$$C_n:=\{m=x^{\alpha}y^{\beta}\notin E\mbox{ such that }
d(m)\leq n\}
$$ Each $m$ of $C_n$ is in one sector $\Delta_i$. Put 
$f_m=P(c_i).\frac{m}{c_i}$.
A division of $f_m$ by the elements $f_{m'}$ smaller than $f_m$ 
with respect to the order $>_-$ is:
$$
f_m=\sum_{d(m')=d(m),\ m'>_+m} q_{m'}f_{m'}+g_m, \mbox{ where
  }q_{m'}\in k.
$$
The remainder $g_m$ writes down 
$$
g_m=m+\sum_{m_i \in E,\ d(m_i)=d(m)} \mu_{m,m_i} m_i.
$$
By construction of $Hilb(\pla)$, if we have chosen 
$n$ big enough, $i\circ e$ is given in Pl\"ucker coordinates 
by $\bigwedge_{m \in C_n}g_m$. The coordinate corresponding to 
the term $\bigwedge_{m \in C_n}m$ equals 1 and $i\circ e$ is then a
morphism from $\s k[X_{c,m}]$ to an affine space $\s k[Y_1,\dots,Y_p]$. 
If $(c,m)$ is a significant cleft couple, the 
coordinate corresponding to the factor
$m \bigwedge (\wedge_{m' \in C_n,m'\neq c}m')$ is $\mu_{c,m}=X_{c,m}+R_{c,m}$,
where $R_{c,m}$ is a polynomial in the variables $X_{c',m'}$, with
$c'<_-c$ or ($c'=c$ and $m'<_+m$). This means that we can order the variables
$X_{c,m}$  and the variables $Y_1,\dots,Y_p$
such that the morphism
$$
i\circ e:\s k[X_{c,m}]=:\s k[X_1,\dots,X_q] \fd \s k[Y_1, \dots,Y_p]
$$
is given by $Y_1=X_1$ and for $1<i\leq q$, $Y_i=X_i+R_i$ where $R_i$ is a 
polynomial in the variables $X_j$, $j<i$. This shows that $i\circ e$
is an embedding. 
\findem

\subsection{The general case}
\label{ssec:exponentielle cas general}
In the previous section, we have defined 
a map $e:T_{Z(E)}^{ab+}\fd \HH_{ab}$ which induces
an isomorphism
$e:T_{Z(E)}^{ab+}\fd C(E)$. In this section, we extend the 
map $e$ to a map  $T_{Z(E)}^+ \fd \HH$, which induces an isomorphism
between $T_{Z(E)}^+$ and a Bialynicki-Birula cell $C(E)$ of $\HH$ 
and we explain the link between this description of $C(E)$ and 
the one of \cite{ellingsrud-stromme88:cartes_sur_le_schema_de_Hilbert}.
Note that this section is not useful for the proof of
theorem \ref{thm principal d'irreductibilite}.
\nl
Choose $a<0$ and  $b<0$ relatively prime. The one-dimensional torus
$k^*=T_{ab}$  acts on $\HH$ with the monomial subschemes as the only
invariants subschemes by an argument similar to \ref{si ab>0 hilbert equiv=ens de points},
so we can apply Bialynicki-Birula to define cells
from this action. To be concrete,  the cell $C(E)$ parametrizes the subschemes $X$
such that $lim_{t \fd 0}(t^{-b},t^a).X=Z(E)$.
\nl
Choose a variable $X_{c,m}$ for each cleft couple with positive
half-direction. Define as in the previous section polynomials
$P(c_i)\in k[X_{c,m}][x,y]$ for each
cleft $c_i$, $Q(c_i,m)\in k[X_{c,m}][x,y]$ for  each significant
$(c_i,m)$. Explicitly, with the notations of the previous section 
$$P(c_n)=c_n$$
$$Q(c_i,m)=P(c_{i+k}).\f{m}{c_{i+k}}$$
and for $i<n$, 
$$P(c_i)=P(c_{i+1}).\f{c_i}{c_{i+1}}+\sum_{m \mbox{ s.t. $(c_i,m)$ is
significant}} X_{c_i,m}Q(c_i,m).$$
\begin{prop}
\label{prop:platitude_de_la_famille-cas general}
The ideal $I_+$ generated by the $P(c_i)$ defines a flat family 
$Z_+$ of constant length over
$\s k[X_{c,m}]$. 
\end{prop}
\demo: this is exactly the same proof as in
\ref{prop:platitude_de_la_famille} except for the claim of flatness
which is not valid since the family is not any longer
supported by the origin. To conclude that $Z_+$ is flat over 
$\s k[X_{c,m}]$ knowing
that the fibers are of constant length, we need to check that if $p$
is the generic point of $\s k[X_{c,m}]$, if $G_p$ is the support
of the fiber $(Z_+)_p$, then the closure 
of $G_p$ in $\s k[X_{c,m}]\x \plp$ is a subscheme of $\s
k[X_{c,m}]\x \pla$. 
\\
The ideal of $G_p$ contains an ideal $(x,y^l+a_{l-1}y^{l-1}+\dots a_0 y^0)$
where $a_i$ is a polynomial in the indeterminates $X_{c,m}$. The
closure $\overline{G_p}$ of $G_p$ in $\s k[X_{c,m}]\x \plp$ is included in the
subscheme defined by the homogeneous ideal $(x,y^l+a_{l-1}y^{l-1}h+\dots
a_0 y^0h^l)$. Thus $\overline{G_p}$ doesn't meet the relative line
$h=0$. 
\findem

\begin{thm}
The morphism $e:T_{Z(E)}^{ab+}\fd \HH$ defined by $I_+$
induces an isomorphism $T_{Z(E)}^{ab+}\simeq C(E)$.
\label{prop:e_est_un_plongement cas general} 
\end{thm}
\demo: the proof is essentially identical to the proof of 
\ref{prop:e_est_un_plongement}, with few changes left to the reader. \findem

\begin{rem} \label{rem:e est un plongement avec tangent negatif}
  We can clearly by the same method parametrize the cells of the action of a
  torus $T_{ab},a>0, b>0$ using the negative tangent space.
\end{rem}

\noindent
We now describe the link between  this last theorem and theorem 2 of 
\cite{ellingsrud-stromme88:cartes_sur_le_schema_de_Hilbert}. 
In the remainder of this section, we assume familiarity with 
the notations and the results 
of \cite{ellingsrud-stromme88:cartes_sur_le_schema_de_Hilbert}. 
\\
In theorem 2 of 
\cite{ellingsrud-stromme88:cartes_sur_le_schema_de_Hilbert},
Ellingsrud and Str\o mme consider the Bialynicki-Birula stratification 
corresponding to an action of a one parameter subgroup
$\psi:\mathbb{G}\fd T$. They
construct a set of variables $s_{ij}$. The one dimensional space
generated by $s_{ij}$ is provided with an action of the torus 
$\psi(\mathbb{G})$ of 
weight $n_{ij}$. Denoting by $S(\psi)$ the polynomial ring 
whose indeterminates are the variables $s_{ij}$ of negative weight,
they construct a family $Z\inc Spec \ S(\psi)\x \pla$, 
flat on $Spec \ S(\psi)$. 
\\
In case $\psi(\mathbb{G})$ coincides with our torus $T_{ab}$, we
can make the link between the two theorems. 
\\
In this particular case where $\psi(\mathbb{G})=T_{ab}$ with $a<0$
and $b<0$, the variables $s_{ij}$ with negative 
weight are the variables $s_{ij}$ with $(i,j)\in \Delta_1$ with 
their notations. 
We make a bijection $\phi$ between the set $R_+$ of significant
positive cleft couples and $\Delta_1$.
If $(c,m)\in R_+$ , let $l>0$ be the smallest integer
such that $x^l.m \notin E$. Using the numeration of the monomials
of \cite{ellingsrud-stromme88:cartes_sur_le_schema_de_Hilbert},
we have $cx^l=d_j$ and $mx^l=d_i$. We let $\phi(c,m)=(i,j)$. 
The following proposition whose proof is omitted says that, up to
sign, $\phi$ identifies 
the family $Z_+$ giving theorem 
\ref{prop:e_est_un_plongement cas general} and the family $Z$ giving theorem 2
of Ellingsrud and Str\o mme. 

\begin{prop}
  There exists an application $\epsilon:R_+\fd \{1,-1\}$ such that,
  if $\eta: \s S(\psi) \fd \s k[X_{c,m}]$ is the morphism defined by
  $X_{c,m}\mapsto \epsilon(c,m)s_{\phi(c,m)}$, then $Z=Z_+\x _{\s
  k[X_{c,m}]}  \s S(\psi)$.
\end{prop}

\begin{rem}
  It seems that the theorem of Ellingsrud and Str\o mme is more general
  than ours, since they can take any one parameter subgroup of $k^*\x
  k^*$ whereas we only consider the sub-tori $T_{ab}$ with
  $ab>0$. However, there is a gap in theorem 2 in the extra cases
  since the family they produce is not always of constant length. 
  For instance, take the one parameter subgroup $t\mapsto (t,t)$, 
  and, following the procedure of
  \cite{ellingsrud-stromme88:cartes_sur_le_schema_de_Hilbert},
  construct a family whose fiber over the origin is the subscheme
  $Z(E)$ of
  length 6, where
  $I^E=(y^3,xy,x^4)$. The fiber of $Z$
  over the point with all coordinates $s_{ij}=0$ but $s_{03}=1$,
  $s_{74}=1$ is of length 7 since a Gr\"obner basis with respect
  to the homogeneous order with $x>y$ is $xy^2+y^3,x^2y+xy^2,x^3+x^2y-xy-y^2,y^4-y^3$.
  In particular, one may see theorem \ref{prop:e_est_un_plongement cas
  general} and remark \ref{rem:e est un plongement avec tangent negatif} 
  as sufficient conditions under which the theorem of Ellingsrud and Str\o mme
  applies. 
  \\
  Let us explain when these conditions apply.
  There are 3 types of cells for $Hilb^l(\plp)$ in 
  \cite{ellingsrud-stromme88:cartes_sur_le_schema_de_Hilbert},
  described in theorems 3,4,5.
  We note that the sufficient conditions given in the present  paper are exactly the 
  conditions needed in theorem 5. One can show that theorem 3 is
  correct too (this concerns the  cells of $\HH$ which
  are the most important cells). The above counterexample concerns
  theorem 4.
\end{rem}

\section{Connecting the strata.}
\label{sec:drawing curves on the tangent space}
Let $a>0$ and $b<0$. The monomials 
of $k[x,y]$ can be ordered in 
an infinite sequence $m_0<m_1<\dots$ with respect to the order $>$ 
defined by $m<m'$ if $(d(m),d_y(m))<(d(m'),d_y(m'))$ for the
lexicographic order on $\NN^2$. For a staircase $E$, 
let $S_E$ be the function
from $\NN$ to $\NN$ defined by $S_E(k)=$number of monomials in  
$E$ smaller or equal to $m_k$.
\begin{defi}
  We define a partial relation $>$ on staircases by: $E>F$ if
  $\forall k$, $S_E(k)>S_F(k)$. 
\end{defi}
\noindent
The goal of this section is to show the following proposition:

\begin{prop} \label{prop:description de la degenerescence}
Let $E$ be a staircase such that
$T_{Z(E)}^{ab+}\neq 0$.
Then there exists a staircase $F$  
verifying $F<E$ and 
$\overline{C(E)} \cap C(F) \neq \emptyset$.   
\end{prop}
\demo: choose a point $p \in \s k[X_{c,m}]$ different from the
origin $O$. There exists a
monomial $m$ which is in $I^E(p)$ but not in $I_+(p)$ ( otherwise 
$I^E(p) \inc I_+(p)$, and even $I^E(p)=I_+(p)$ since the families defined by 
$I^E(p)$ and $I_+(p)$ have the same relative degree, which is impossible
since $e(O)\neq e(p)$ by \ref{prop:e_est_un_plongement} ). 
Take such an $m$ minimal.
The element $g_m(p) \in I_+(p)$ defined in theorem \ref{prop:e_est_un_plongement}
has its terms in $E$ except the initial term $m$.
There is at least one term different from $m$ since $m\notin I_+(p)$ by definition.
Consider the subscheme $X(p)$ 
defined by the ideal $I_+(p)$. If the torus  
$k^*$ acts on $\s k[x,y]$ by $t.x=tx$, $t.y=y$, then the scheme 
$X(\infty)=\lim_{t \fd \infty}t.X(p)$ is not in $C(E)$ since its ideal
contains the monomial $\lim_{t \fd \infty} t.g_m(p)\in E$.  If $F$ is
the staircase of $X(\infty)$, then $\overline{C(E)}\cap C(F)\neq
\emptyset$. This non vacuity implies $E>F$
(\cite{evain00:cnincidence_cellules_schubert}
or \cite{yameogo94:cn_incidence_compositio}). 
\findem

\section{Uniqueness of the final stratum and conclusion of the proof}
\label{sec:unicity of the final point}
The connectedness of $\HH_{ab}(H)$ has been proved in the first section under the 
condition $ab\leq 0$. In the remaining cases, one can suppose
$a>0$ and $b<0$ and we have stratified $\HH_{ab}(H)$ by 
affine spaces $C(E)$,
ordered the staircases parametrizing the strata 
and proven in the previous section
that $\overline{C(E)} \cap C(F) \neq \emptyset$ for some staircase $F<E$
provided that 
$T_{Z(E)}^{ab+}\neq 0$. We can go  down from stratum to stratum while the positive
tangent space is non trivial. The process stops 
after a finite
number of steps since the number of possible staircases is finite.
In fact, the stratum $C(E_m)$ in which the process stops doesn't
depend on the intermediate degenerations. The staircase $E_m$
parametrizing this final stratum is characterized by the following theorem. 
\\
Recall that a staircase $E$ is compatible with $H=(\dots,h_{-1},h_0,h_1,\dots)$ if
$C(E)\cap H_{ab}(H)\neq \emptyset$ or equivalently if $E$ has $h_i$
elements in degree $i$.
\begin{thm} \label{thm:description escalier final}
  Let $S$ be the set of staircases compatible with $H$. If $S\neq \emptyset$, then
  there exists in $S$  a staircase $E_m$ such that: $\forall E\in S, \
  E_m\leq E$.
\end{thm}
\noindent
To conclude the proof of the connectedness, 
it clearly suffices to prove this last theorem
and to show that the path through the strata 
always stops in that special stratum $C(E_m)$.
\\
Since a minimal staircase $E$ verifies $T_{Z(E)}^{ab+}=0$
by the previous section, and since we stop when $T_{Z(E)}^{ab+}=0$,
we just have to show that there is a unique staircase $E \in S$ such that
$T_{Z(E)}^{ab+}=0$, which we prove by induction on the cardinal of $E$.
\nl
The case $card(E)=1$ is trivial. 
A general $E$ is the ``vertical collision'' of 
the bottom row $L$ and of a residual staircase $U$
determined by its  monomials: 
$E_{R}:=\{m \mbox{ such that } ym \in E\}$. 
It is in fact sufficient to show that the bottom
row $L$ of $E$ is completely determined. 
Indeed,  suppose that $L$ is known. The
positive tangent space at $Z(U)$ injects in the positive tangent space at 
$Z(E)$ by sending $\sum \lambda_{c,m} (c,m)$ to $\sum \lambda_{c,m}
(yc,ym)$, so $T^{ab+}_{Z(U)}=T^{ab+}_{Z(E)}=0$. 
It follows that $U$ is completely determined by induction (since its
Hilbert function is determined by that of $E$ and $L$).
\nl
The next two lemmas complete the proof. They explain 
that the bottom row $L$ of a staircase $E$ with trivial positive tangent
space and compatible with $H$ is determined: the maximal integer
$k$ such that $x^k\in E$ is the maximum integer such that
$H(kb-a)<H(kb)$. More precisely, lemma \ref{longueur_mini_ligne} 
shows that the bottom 
row contains all the monomials $x^k$ such that $H(kb-a)<H(kb)$
and lemma \ref{longueur_maxi_ligne} 
shows that it cannot contain a bigger monomial. 
\\

\begin{lm}\label{longueur_maxi_ligne}
  Let $E$ be a staircase compatible with $H$. Let $x^k$ be the maximal
  power of $x$ such that $x^k\in E$ and let $\delta=kb$ be its degree. If
  $H(\delta)\leq H(\delta-a)$, then $T^{ab+}_{Z(E)}\neq \{0\}$.
\end{lm}
\demo: E contains a monomial $m$ of degree $\delta$ 
such that $ym\notin E$. Otherwise, $\phi:m\mapsto ym$ would be an
injective application from the set $E(\delta-a):=\{m \in E, \ d(m)=b\}$
to $E(\delta)$. The map $\phi$ is not surjective since $x^k \notin
Im(\phi)$. Thus, $card(E(\delta))=H(\delta)>H(\delta-a)=card(E(\delta-a))$, contradicting
the hypothesis. \\
Now choose $l\in \NN$ maximal such that $c=\f{ym}{x^l} \in k[x,y]-E$. 
By construction, $(c,x^{k-l})$ is a cleft couple, which shows 
$T_{Z(E)}^{ab+} \neq 0$. 
\findem

\begin{lm}\label{longueur_mini_ligne}
  Let $\delta=kb$ be an integer. If $H(\delta-a)<H(\delta)$, 
  then for any staircase $E$
  compatible with $H$, the  element $x^k$  of
  degree $\delta$ is in the bottom row of $E$.
\end{lm}
\demo:if $x^k$ were not in $E$,
then to each element $x^{\alpha}y^{\beta} \in E$
of $E$ of degree $\delta$
would be associated the element $x^{\alpha-1}y^{\beta} \in E$ of degree $\delta-a$. A
count of these elements would then show
$H(\delta-a)\geq H(\delta)$.\findem
\nl

\begin{rem}
The existence of a minimal staircase with respect to the partial order
on staircases 
in theorem \ref{thm:description escalier final} is a purely combinatorial
statement which has required an algebro-geometric proof. I have not found a
combinatorial proof simpler than the given argument.
\end{rem}


\begin{thebibliography}{10}

\bibitem{bialynicki-birula73:theorems_on_actions_of_groups_annals}
A~Bialynicki-Birula.
\newblock Some theorems on actions of algebraic groups.
\newblock {\em Ann. of Math.}, (98):480--497, 1973.

\bibitem{briancon77:description_de_Hilb_a2}
J~Brian\c{c}on.
\newblock Description de ${H}ilb^n\mathbb{C}\{X,y\}$.
\newblock {\em Invent .Math.}, (41), 1977.

\bibitem{brion97:_equivariant_chow_groups}
M~Brion.
\newblock Equivariant chow groups for torus actions.
\newblock {\em Transformations Groups}, (2):225--267, 1997.

\bibitem{ellingsrud-stromme87:chow_group_of_hilbert_schemes}
Geir Ellingsrud and Stein Arild~Str\o mme.
\newblock On the homology of the {H}ilbert scheme of points in the plane.
\newblock {\em Invent. math.}, (87):343--352, 1987.

\bibitem{ellingsrud-stromme88:cartes_sur_le_schema_de_Hilbert}
Geir Ellingsrud and Stein Arild~Str\o mme.
\newblock On a cell decomposition of the {H}ilbert scheme of points in the
  plane.
\newblock {\em Invent. Math.}, (91):365--370, 1988.

\bibitem{evain00:cnincidence_cellules_schubert}
L~Evain.
\newblock Incidence relations among the {S}chubert cells of equivariant
  {H}ilbert schemes.
\newblock {\em Math. Z.}, to appear.

\bibitem{fogarty68:lissite_du_hilbert}
J~Fogarty.
\newblock Algebraic families on an algebraic surface.
\newblock {\em Am. J. Math.}, (10):511--521, 1968.

\bibitem{fulton93:varietes-toriques}
W~Fulton.
\newblock {\em Introduction to toric varieties}, volume 131.
\newblock Princeton University Press, 1993.

\bibitem{hartshorne66:connexite_du_schema_de_hilbert}
R~Hartshorne.
\newblock Connectedness of the {H}ilbert scheme.
\newblock {\em Publ. Math. I.H.E.S}, (29):5--48, 1966.

\bibitem{iarrobino77:_punctual_hilbert_schemes_AMS}
A~Iarrobino.
\newblock Punctual {H}ilbert schemes.
\newblock {\em Mem. Am. Math. Soc.}, (188), 1977.

\bibitem{iversen72:lissite-de-partie-invariante-sous-l'action-d'un-tore}
B.~Iversen.
\newblock A fixed point formula for action of tori on algebraic varieties.
\newblock {\em Inventiones math.}, (16):229--236, 1972.

\bibitem{lehn99:_chern_classes_of_tautological_sheaves_on_Hilbert_schemes}
M~Lehn.
\newblock Chern classes of tautological sheaves on hilbert schemes of points on
  surfaces.
\newblock {\em Invent. math.}, (136):157--207, 1999.

\bibitem{nakajima97:_heisenberg_et_Hilbert_schemes}
H~Nakajima.
\newblock Heisenberg algebra and hilbert schemes of points on projective
  surfaces.
\newblock {\em Ann. of Math.}, 2(145):379--388, 1997.

\bibitem{yameogo94:cn_incidence_compositio}
J~Yameogo.
\newblock D\'ecomposition cellulaire de vari\'et\'es param\'etrant des id\'eaux
  homog\`enes de $\mathbb{C}[[x,y]]$.
\newblock {\em Compositio Math.}, 1(90):81--98, 1994.

\end{thebibliography}

\end{document}